\documentclass[12pt]{article}
\usepackage{amsmath,amssymb,amsfonts,theorem}
\newtheorem{thm}{Theorem}[section]
\newtheorem{lem}[thm]{Lemma}
\newtheorem{prop}[thm]{Proposition}
\newtheorem{cor}[thm]{Corollary}
 {
 \theorembodyfont{\rmfamily}
 
 \newtheorem{ex}[thm]{Example}
 }
\newenvironment{pf}{\paragraph{Proof}}{\par\medskip}
\newcommand{\qed}{\ifhmode\unskip\nobreak\fi\quad\ensuremath\diamond}
\newcommand{\iso}{\cong}
\newcommand{\into}{\hookrightarrow}
\newcommand{\tensor}{\otimes}
\newcommand{\C}{\mathbb C}
\newcommand{\Z}{\mathbb Z}

\newcommand{\pp}{\mathbb P}
\DeclareMathOperator{\Aut}{Aut}
\DeclareMathOperator{\Pic}{Pic}
\newcommand{\Oh}{\mathcal O}
\newcommand{\epsi}{\varepsilon}

\newcommand{\al}{\alpha}

\newcommand{\ga}{\gamma}
\newcommand{\Ga}{\Gamma}
\newcommand{\De}{\Delta}
\newcommand{\Si}{\Sigma}

\newcommand{\fie}{\varphi}
\newcommand{\Phat}{\widehat{\pp}}
\newcommand{\Dbar}{\overline D}
\newcommand{\fbar}{\overline f}

\newcommand{\inv}{^{-1}}
\renewcommand{\labelenumi}{(\roman{enumi})}
\numberwithin{equation}{section}

\begin{document}
\title{The bicanonical map of surfaces \\ with $p_g=0$ and $K^2\ge7$
\thanks{2000 Mathematics Subject Classification: 14J29, 14E05}}
\author{Margarida Mendes Lopes \and Rita Pardini}
\date{}

\maketitle
\begin{abstract} A minimal surface of general type with
$p_g(S)=0$
 satisfies $1\le
K^2\le 9$ and  it is known that the image of the
bicanonical map $\fie$ is a surface for $K_S^2\geq 2$, whilst for $K^2_S\geq
5$,  the bicanonical map is always a morphism. In this paper it is shown  that 
 $\fie$ is birational if  $K_S^2=9$ and that
 the degree of $\fie$ is at most $2$  if $K_S^2=7$ or
$K_S^2=8$.  \par By presenting two examples of surfaces $S$ with $K_S^2=7$ and $8$
and bicanonical map of degree 2, it is also shown that this result is sharp. The example with $K_S^2=8$ is, to our knowledge, a
 new example of a surface
of general type with $p_g=0$.  \par The
degree of
$\fie$ is also calculated  for two other
known surfaces of general type with $p_g=0$, $K_S^2=8$. In both cases the
bicanonical map turns out to be birational.
\end {abstract}
\section{Introduction}

Many examples of complex surfaces of general type with $p_g=q=0$ are known,
but a detailed classification is still lacking, despite much progress in the
theory of algebraic surfaces. Surfaces of general type are often studied
using properties of their canonical curves. If a surface has $p_g=0$, then
there are of course no such curves, and it seems natural to look instead at
the bicanonical system, which is not empty.

Minimal surfaces $S$ of general type with $p_g(S)=0$ satisfy $1\le
K_S^2\le9$. By a result of Xiao Gang \cite{xiaocan}, the image of the
bicanonical map is a surface for $K_S^2\ge2$, while for $K^2_S\ge5$, by
Reider's theorem \cite{reider}, the bicanonical map is always a morphism.
Xiao Gang \cite{xiao2} showed that the degree of the bicanonical map is
$\le2$ for surfaces of general type, with a limited number of possible
exceptions. Surfaces with $p_g=0$ are among the exceptional cases, and
\cite{xiao2} gives practically no information on the possible degrees
of their bicanonical maps.

The first author \cite{marg} showed that if $K_S^2\ge3$ and the
bicanonical map is a morphism, its degree is $\le4$. There are examples
due to Burniat \cite{bu} (see also \cite{peters}) with $3\le K^2\le 6$
having bicanonical map of degree $4$; indeed, \cite{MP} gives a precise
description of the surfaces with $K_S^2=6$, $p_g(S)=0$ and bicanonical
map of degree 4. Here we refine the result of \cite{marg} by proving the
following results:

 \begin{thm}\label{main}
Let $S$ be a minimal surface of general type defined over $\C$ with
$p_g(S)=0$, and $\fie\colon S\to \Si\subset\pp^{K^2_S}$ its bicanonical
map, with image $\Si$.
 \begin{enumerate}
\item If $K^2_S=9$ then $\fie$ is birational;
 \item if $K_S^2=7,8$ then $\fie$ has degree $\le 2$.
\end{enumerate}
\end{thm}

\begin{prop}\label{sharp}
There exist minimal surfaces of general type with $p_g(S)=0$, $K_S^2=7,8$
and bicanonical map of degree $2$.
 \end{prop}

We prove this by giving two examples of surfaces $S$ with $K_S^2=7$ and $8$
and bicanonical map of degree 2. The example with $K_S^2=7$ is due to Inoue
\cite[Remark~6]{inoue}, who constructed it as a quotient of a complete
intersection in the product of four elliptic curves by a free action of
$\Z_2^5$. Here we give an alternative description as a $\Z_2^2$-cover of a
singular rational surface that allows us to describe the bicanonical map and
compute its degree. The example with $K_S^2=8$ is obtained by applying a
construction of Beauville (\cite[p.~123, Exercise~4]{bv}, and cf.\
\cite{dolg}). To the best of our knowledge, it is a new example of a surface
of general type with $p_g=0$.

Only a few examples of surfaces with $p_g=0$, $K_S^2=8$ are known. It is a
nontrivial exercise to compute the degree of the bicanonical map for an
explicit surface, and for interest, we include the computation for two other
known surfaces of general type with $p_g=0$, $K_S^2=8$, for both of which the
bicanonical map turns out to be birational.

The paper is organized as follows: Section~\ref{covers} recalls some facts
on irregular double covers, one of the main ingredients of the proof of the
theorem; in Section~\ref{teorema} we prove the main result. For $K^2_S=7,8$
the proof consists of using the methods of Section~\ref{covers} to exclude
the possibility that the bicanonical map has degree 4; for $K^2_S=9$ the
result is obtained by combining Reider's theorem with an analysis of the
Picard group of $S$. In the final Section~\ref{exa} we present the two
examples to prove Proposition~\ref{sharp}, and we also compute the degree of
the bicanonical map of two other surfaces with $p_g=0$ and $K^2_S=8$.

\paragraph{Notations and conventions} We work over $\C$; all varieties are
assumed to be compact and algebraic. We do not distinguish between line
bundles and divisors on a smooth variety, and use additive and multiplicative
notation interchangeably. Linear equivalence is denoted by $\equiv$ and
numerical equivalence by $\sim$. The remaining notation is standard in
algebraic geometry.

\section{Irregular double covers and fibrations}\label{covers}

We describe here the key facts used in some proofs in this paper.

Let $S$ be a smooth complex surface, $D\subset S$ a curve (possibly empty)
with at worst ordinary double points, and $M$ a line bundle on $S$ with
$2M\equiv D$. It is well known that there exists a normal surface $Y$ and
a finite degree $2$ map $\pi\colon Y\to S$ branched over $D$ such that
$\pi_*\Oh_Y=\Oh_S\oplus M\inv$. The singularities of $Y$ are $A_1$ points
and occur precisely above the singular points of $D$; thus it makes sense
to speak of the canonical divisor, the geometric genus, the irregularity
and the Albanese map of $Y$. We refer to $Y$ as the {\em double cover
defined by the relation $2M\equiv D$}. The invariants of $Y$ are:
 \begin{eqnarray}\label{formule}
 K_Y^2& =&2(K_S+M)^2,\nonumber \\
 \chi(\Oh_Y)&= &2\chi(\Oh_S)+\frac{1}{2}M(K_S+M), \\
 p_g(Y)&=&p_g(S)+h^0(S,K_S+M).\nonumber
 \end{eqnarray}
If $p_g(S)=q(S)=0$, the existence of a double cover $\pi\colon Y\to S$ with
$q(Y)>0$ forces the existence of a fibration $f\colon S\to\pp^1$ such that
$\pi\inv$ of the general fibre of $f$ is disconnected. More precisely we have:

 \begin{prop}[De Franchis]\label{defranchis}
 Let $S$ be a smooth surface with $p_g(S)=q(S)=0$ and $\pi\colon Y\to S$ a
double cover with at most $A_1$ points; if $q(Y)>0$, then
 \begin{enumerate}
 \item the Albanese image of\/ $Y$ is a curve $B$;
 \item let $\al\colon Y\to B$ be the Albanese fibration. Then there exists a
fibration $g\colon S\to\pp^1$ and a degree $2$ map $p\colon B\to\pp^1$ such
that $p\circ\al=g\circ\pi$.
 \end{enumerate}
\end{prop}

The possibility of existence of such a double cover often leads to a contradiction, using the following
result:

\begin{cor}\label{genus2} Let $S$ be a minimal surface of general type
with $p_g(S)=q(S)=0$ and $K_S^2\ge3$, and $\pi\colon Y\to S$ a double
cover with at most $A_1$ points. Then $K_Y^2\ge16(q(Y)-1)$.
\end{cor}

Proposition~\ref{defranchis} is an old result of De Franchis \cite{defra},
explained and generalized in several ways by Catanese and Ciliberto
\cite{cetraro}. Proposition~\ref{defranchis} and Corollary~\ref{genus2}
are both stated and proved in \cite{MP} for smooth $Y$, but the proof
extends verbatim to the case of $A_1$ points.

\section{Proof of Theorem~\ref{main}}\label{teorema}
Under the assumptions of Theorem~\ref{main}, the image of the bicanonical
map is a surface by \cite{xiaocan}, and the bicanonical map is a morphism
by Reider's theorem \cite{reider}. Moreover, since $4K^2_S=\deg\fie\deg\Si$
and $\Si$ is a nondegenerate surface in $\pp^{K_S^2}$, the possible values
of $\deg\fie$ are $1,2,4$ for $K^2_S=7,8$ and $1,2,3,4$ for $K^2_S=9$.

We prove the theorem by analysing separately the cases $K^2_S=7,8,9$. In
each case we argue by contradiction.

\subsection{The case $K^2_S=7$}\label{K7}
By the above remark, it is enough to show that $\deg\fie=4$ does not occur.
Assume that $\fie$ has degree 4. The bicanonical image $\Si$ is a linearly
normal surface of degree $7$ in $\pp^7$ and its nonsingular model has
$p_g=q=0$. By \cite[Theorem~8]{nagata}, $\Si$ is the image of the blowup
$\Phat$ of $\pp^2$ at two points $P_1,P_2$ under its anticanonical map
$f\colon \Phat\into\pp^7$. If $P_1\ne P_2$, then $f$ is an embedding,
while if $P_2$ is infinitely near to $P_1$ (say) then $\Si$ has an $A_1$
singularity. In either case, the hyperplane section of $\Si$ can be written as
$H\equiv 2l+l_0$, where $l$ is the image on $\Si$
of a general line of $\pp^2$ and $l_0$ is the image on $\Si$ of the strict
transform of the line through $P_1$ and $P_2$. Notice that $l_0$ is
contained in the smooth part of $\Si$. Thus we have $2K_S\equiv 2L+L_0$,
where $L=\fie^*l$ and $L_0=\fie^*l_0$.

\begin{lem}\label{pullback}
$L_0$ satisfies one of the following possibilities:
 \begin{enumerate}
 \item there exists an effective divisor $D$ on $S$ such that $L_0=2D$; or
 \item $L_0$ is a smooth rational curve with $L_0^2=-4$; or
 \item there exist smooth rational curves $A$ and $B$ with $A^2=B^2=-3$,
$AB=1$, and $L_0=A+B$.
 \end{enumerate}
\end{lem}
\begin{pf} Remark first that $K_SL_0=2$, $L^2_0=-4$, and $L_0=2(K_S-L)$ is
divisible by 2 in $\Pic S$. Let $\theta$ be a $-2$-curve of $S$; then
$\theta$ is contracted by $\fie$ and thus $L\theta=L_0\theta=0$. Since $L$
and $L_0$ are independent elements of the 3-dimensional space $H^{1,1}(S)$,
$S$ contains at most one $-2$-curve. We write $L_0=C+a\theta$, where $C$ is
the strict transform of $L_0$, $\theta$ is a $-2$-curve and $a\ge0$ (we set
$a=0$ if $S$ has no $-2$-curve). The equalities $\theta L_0=0$ and
$L^2_0=-4$ imply
 \begin{equation}\label{conti}
 \theta C=2a, \quad\hbox{and}\quad C^2=-4-2a^2.
 \end{equation}
If $C$ is irreducible, then $K_SC=2$ implies $C^2\ge-4$ and thus $a=0$ and
case (ii) holds. If $C$ is reducible, then $C=A+B$, with $A$ and $B$
irreducible curves such that $K_SA=K_SB=1$. If $A=B$, then
$AL_0=2A^2+a\theta A=2A^2+a^2$ is even, because $L_0$ is divisible by 2, and
thus $a$ is even and we are in case (i). If $A\ne B$, then $AB\ge0$ and
$A^2,B^2\ge-3$; by parity considerations and (\ref{conti}) we get
$A^2=B^2=-3$ and either $AB=1$, $a=0$ or $AB=0$, $a=1$. The first case
corresponds to (iii), while the second does not occur. In fact the
intersection matrix of $A$, $B$, $\theta$ would be negative definite,
contradicting the index theorem, since $h^{1,1}(S)=3$. \qed \end{pf}

In cases (ii) or (iii) of Lemma~\ref{pullback}, let $\pi\colon Y\to S$ be
the double cover given by $2(K_S-L)\equiv L_0$; then
the formulas (\ref{formule}) give
$\chi(Y)=2$ and $K_Y^2=16$. Since the bicanonical map $\fie$ maps $L$ onto
a twisted cubic, $h^0(S,\Oh_S(2K_S-L))=4$ and thus
$p_g(Y)=p_g(S)+h^0(S,\Oh_S(2K_S-L))=4$; we thus obtain $q(Y)=3$,
contradicting Corollary~\ref{genus2}.

In case (i) of Lemma~\ref{pullback}, consider the \'etale double cover
$\pi\colon Y\to S$ given by $2(K_S-L-D)\equiv 0$; arguing as above, we get
that the invariants of $Y$ are
 \[
 K^2_Y=14,\quad \chi(\Oh_Y)=2,\quad p_g(Y)=p_g(S)+h^0(S,\Oh_S(2K_S-L-D))=3,
 \]
so that $q(Y)=2$ and we again obtain a contradiction to
Corollary~\ref{genus2}.

Hence $\deg\fie\ne4$ and we have proved Theorem~\ref{main} in case
$K^2_S=7$.

\subsection{The case $K^2_S=8$}\label{K8}

As in case $K_S^2=7$, it is enough to show that $\deg\fie=4$ does not occur.
If $\fie$ has degree 4, then the bicanonical image $\Si$ is a linearly normal
surface of degree 8 in $\pp^8$ whose nonsingular model has $p_g=q=0$. By
\cite[Theorem~8]{nagata}, $\Si$ is either the Veronese embedding in $\pp^8$
of a quadric $Q\subset\pp^3$ or the image of the blowup $\Phat$ of $\pp^2$ at
a point $P$ under its anticanonical map $f\colon \Phat\into\pp^8$.

In the first case $2K_S\equiv 2A$, where $A$ is the hyperplane section of
$Q$. Then $\eta=K_S-A$ is a nontrivial 2-torsion element in $\Pic S$, since
$p_g(S)=0$. The \'etale double cover $\pi\colon Y\to S$ given by
$2\eta\equiv 0$ has invariants $\chi (Y)=2$, $K_Y^2=16$. Moreover,
$p_g(Y)=p_g(S)+h^0(S,\Oh_S(A))=4$, so that $q(Y)=3$. Since $K_Y^2=16$,
this contradicts Corollary~\ref{genus2}, and therefore $\Si$ is not the
Veronese embedding of a quadric.

If the bicanonical image $\Si$ is the image of $\Phat$ via the map
induced by $|{-}K_{\Phat}|$, then the hyperplane section of $\Si$ can
be written as $H\equiv 2l+l_0$, where $l$ is the image on $\Si$ of a
general line of $\pp^2$ and $l_0$ is the image on $\Si$ of the strict
transform of a general line through $P$. Thus $2K_S\equiv 2L+L_0$, where
$L=\fie^*l$ and $L_0=\fie^*l_0$, and $L_0=\fie^*l_0$ is smooth by
Bertini's theorem. Consider now the double cover $\pi\colon Y\to S$ given
by $2(K_S-L)\equiv L_0$; the formulas (\ref{formule}) give
$\chi(Y)=3$ and $K_Y^2=24$. Since
$p_g(Y)=p_g(S)+h^0(S,\Oh_S(2K_S-L))=
0+h^0(S,\Oh_S(L+L_0))=5$, we get
$q(Y)=3$, contradicting Corollary~\ref{genus2}. Thus $\Si$ is also not
the image of $\Phat$.

Hence $\deg\fie\ne4$ and the proof of Theorem~\ref{main}, (ii) is complete.

\subsection{The case $K^2_S=9$}\label{K9}

If $K^2_S=9$, then by Poincar\'e duality, $H^2(S,\Z)$ is generated up to
torsion by the class of a line bundle $L$ with $L^2=1$; thus every divisor
on $S$ is numerically a multiple of $L$, and in particular $K_S\sim 3L$.

Assume by contradiction that $\fie$ is not birational; then by Reider's
theorem (cf.\ \cite[Theorem~2.1]{bs}), for every pair of points $x_1,
x_2\in S$ with $\fie(x_1)=\fie(x_2)$ there exists an effective divisor
$C$ containing $x_1,x_2$ such that $K_S C-2\le C^2<\frac{1}{2}K_S C<2$.
Since $K_S\sim 3L$, the only possibility is that $C\sim L$. We can assume
that, as $x_1$ and $x_2$ vary, the divisor $C$ varies in an irreducible
system of curves, which is linear by the regularity of $S$. Every curve of
$|C|$ is irreducible, since the class of $C$ generates $H^2(S,\Z)$ up to
torsion, and the general curve of $|C|$ is smooth by Bertini's theorem,
since $C^2=1$. Therefore $|C|$ is a linear pencil of curves of genus~3 with
one base point. For a general $C\in |C|$ we consider the exact sequence:
 \begin{equation}\label{seq}
 0\to \Oh_S(2K_S-C)\to\Oh_S(2K_S)\to\Oh_C(2K_S)\to 0.
 \end{equation}
Since $2K_S-C\sim K_S+2L$, Kodaira vanishing gives $H^1(S,\Oh_S(2K_S-C))=0$,
and the map 
$H^0(S,\Oh_S(2K_S))\to H^0(C,\Oh_C(2K_S))$ induced by the sequence (\ref{seq}) 
is surjective. So the map $f\colon C\to\pp^3$
given by $|\Oh_C(2K_S)|$ is not birational; it follows that $f$ maps $C$
two-to-one onto a twisted cubic, and thus $C$ is hyperelliptic. If we
denote by $\De$ the $g^1_2$ of $C$, then $2K_S|_C\equiv 3\De$ and also, by
the adjunction formula, $K_S+C|_C\equiv 2\De$. So
$\eta\equiv 4K_S-(3K_S+3C)\equiv K_S-3C$ is trivial when restricted to $C$.
Moreover $\eta\sim 0$ and so $\eta$ is a torsion element of $\Pic S$. Since
$p_g(S)=0$, $\eta$ is nonzero. Consider the connected \'etale
cover $\pi\colon Y\to S$ associated to $\eta$. Because $\eta|_C=0$, the
cover $\pi|_{\pi\inv(C)}\colon\pi\inv(C)\to C$ is trivial and thus
$\pi\inv(C)$ is a smooth disconnected curve with each component of
self-intersection 1. This contradicts the Index theorem and we have thus
proved Theorem~\ref {main}, (i). \qed

\section{Examples}\label{exa}

This section calculates the degree of the bicanonical map in 4 interesting
examples, as discussed in the introduction.

\begin{ex}\label{ex!1} Starting from the quadrilateral $P_1P_2P_3P_4$ in
$\pp^2$ of Figure~\ref{fig!quad}, let $P_5$ be the intersection point of
the lines $P_1P_2$ and $P_3P_4$ and $P_6$ the intersection point of $P_1P_4$
and $P_2P_3$. Write $\Si\to\pp^2$ for the blowup of $P_1,\dots,P_6$, and
$e_i$ for the exceptional curves of $\Si$ over $P_i$. Denote by $l$ the
pullback of a line.

Write $S_1,\dots,S_4$ for the strict transforms on $\Si$ of the sides
$P_iP_{i+1}$ of the quadrilateral $P_1P_2P_3P_4$ (we take subscripts
modulo 4); these are the only $-2$-curves of $\Si$. The morphism
$f\colon\Si\to\pp^3$ given by $|{-}K_{\Si}|$ has image a cubic surface
$V\subset \pp^3$, and $f$ is an isomorphism on $\Si\setminus\bigcup S_i$,
and contracts each $S_i$ to an $A_1$ point.

If $A\subset\{P_1,\dots,P_6\}$ consists of 4 points no three of which are
collinear, then the linear system of conics through the points of $A$ gives
rise to a free pencil on $\Si$; we denote by $f_1$ the strict transform of
a general conic through $P_2P_4P_5P_6$, by $f_2$ that of a general conic
through $P_1P_3P_5P_6$ and by $f_3$ that of a general conic through
$P_1P_2P_3P_4$.

Finally, we introduce the ``diagonals'' of the quadrilateral
$P_1P_2P_3P_4$, writing $\De_1,\De_2,\De_3$ for the strict transform
of $P_1P_3$, $P_2P_4$ and $P_5P_6$.
 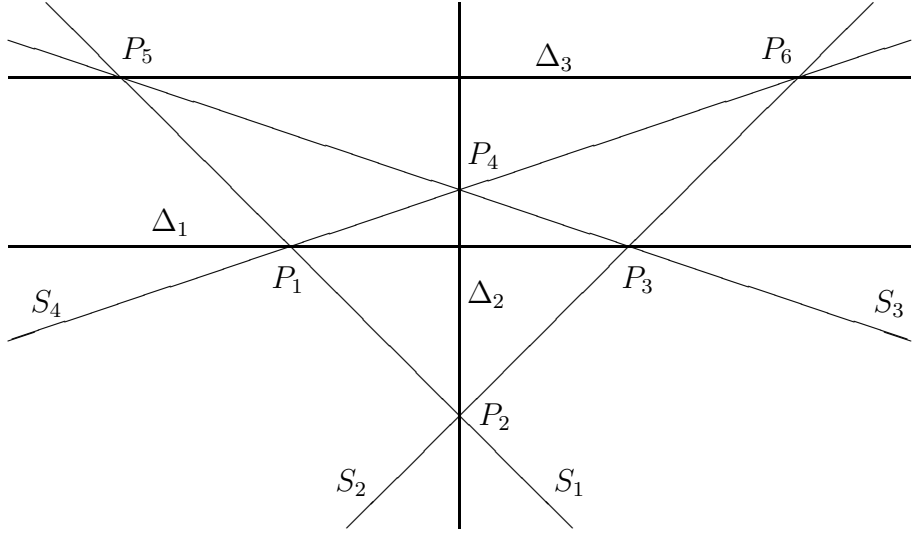
\begin{figure}[ht]
\setlength{\unitlength}{0.5cm}
\centerline{\begin{picture}(24,14)
\thicklines
\put(0,12){\line(1,0){24}}
\put(12,0){\line(0,1){14}}
\put(0,7.5){\line(1,0){24}}
\thinlines
\put(0,13){\line(3,-1){24}}
\put(24,13){\line(-3,-1){24}}
\put(1,14){\line(1,-1){14}}
\put(23,14){\line(-1,-1){14}}
\put(3,12.5){$P_5$}
\put(20,12.5){$P_6$}
\put(7,6.4){$P_1$}
\put(16.3,6.4){$P_3$}
\put(12.5,2.75){$P_2$}
\put(12.2,9.7){$P_4$}
\put(14,12.3){$\De_3$}
\put(12.2,6){$\De_2$}
\put(3.8,7.9){$\De_1$}
\put(0.6,5.75){$S_4$}
\put(23,5.75){$S_3$}
\put(8.7,1){$S_2$}
\put(14.5,1){$S_1$}
\end{picture}}
 \caption{The quadrilateral $P_1P_2P_3P_4$ in $\pp^2$} \label{fig!quad}
 \end{figure}
The divisors we have introduced satisfy the following relations:
 \begin{enumerate}
 \item ${-}K_{\Si}\equiv \De_1+\De_2+\De_3$;
 \item $f_i\equiv \De_{i+1}+\De_{i+2}$ for all $i\in \Z_3$;
 \item $\De_iS_j=0$ for all $i,j$;
\item $\De_if_j=2\delta_{ij}$ for $1\le i,j\le 3$.
 \end{enumerate}
Denote by $\ga_1,\ga_2,\ga_3$ the nonzero elements of $\Ga=\Z_2\times\Z_2$
and by $\chi_i\in\Ga^*$ the nontrivial character orthogonal to $\ga_i$; by
 \cite[Propositions~2.1 and~3.1]{ritaabel}, to
define a smooth $\Ga$-cover $\pi\colon X\to \Si$, we specify:
 \begin{enumerate}
 \renewcommand{\labelenumi}{(\Roman{enumi})}
 \item smooth divisors $D_i$ for $i=1,2,3$ such that $D=D_1+D_2+D_3$ is a
normal crossing divisor,
 \item line bundles $L_1$, $L_2$ satisfying $2L_1\equiv D_2+D_3$,
$2L_2\equiv D_1+D_3$.
 \end{enumerate}

The branch locus of $\pi$ is $D$. More precisely, $D_i$ is the image of the
divisorial part of the fixed locus of $\ga_i$ on $S$. We have
 \[
 \pi_*\Oh_S=\Oh_{\Si}\oplus L_1\inv\oplus L_2\inv\oplus L_3\inv,
 \]
where $L_3=L_1+L_2-D_3$, and $\Ga$ acts on $L_i\inv$ via the character
$\chi_i$.

Here we set:
 \begin{enumerate}
 \renewcommand{\labelenumi}{(\Roman{enumi})}
 \item $D_1=\De_1+f_2+ S_1+S_2$,
 $D_2=\De_2+f_3$,
 $D_3=\De_3+f_1+f_1'+S_3+S_4$;
where $f_1, f_1'\in |f_1|$, $f_2\in |f_2|$, $f_3\in |f_3|$ are
general curves;
\item $L_1=5l-e_1-2e_2-e_3-3e_4-2e_5-2e_6$, and

$L_2=6l-2e_1-2e_2-2e_3-2e_4-3e_5-3e_6$
\end{enumerate}
and we obtain $L_3=4l-2e_1-2e_2-2e_3-e_4-e_5-e_6$.
For $i=1,\dots,4$, the (set theoretic) inverse image of $S_i$ in $X$
is the disjoint union of two $-1$-curves $E_{i1}$, $E_{i2}$; contracting
these $8$ exceptional curves on $X$ and contracting the $S_i$ on $\Si$, we
obtain a smooth $\Z_2^2$-cover $p\colon S\to V$. The map $p$ is branched on
the four singular points of $V$ and on the image $\Dbar$ of $D$, which is
contained in the smooth locus of $V$. The bicanonical divisor $2K_X$ is equal
to $\pi^*(2K_{\Si}+D)=\pi^*(-K_{\Si}+f_1+S_1+S_2+S_3+S_4)=
\pi^*(-K_{\Si}+f_1)+2\sum E_{ij}$, and thus the bicanonical divisor $2K_S$
is equal to $\pi^*(-K_V+\fbar_1)$, where $\fbar_1$ is the image of $f_1$
in $V$. So $2K_S$ is ample, since it is the pullback of an ample line bundle
by a finite map, $S$ is minimal and of general type, and
$K^2_S=\frac{1}{4}4(K_V+\fbar_1)^2=7$.

To compute the geometric genus of $S$, recall that
$p_g(X)=p_g(\Si)+\sum h^0(\Si, K_{\Si}+L_i)$
(cf.\ \cite{quatro} or \cite[Lemma~4.2]{ritaabel}). We have
 \begin{align*}
 K_{\Si}+L_1 &=2l-e_2-2e_4-e_5-e_6, \\
 K_{\Si}+L_2 &=3l-e_1-e_2-e_3-e_4-2e_5-2e_6, \\
 K_{\Si}+L_3 &=l-e_1-e_2-e_3.
 \end{align*}
We show that $h^0(\Si,K_{\Si}+L_2)=0$. Assume
by contradiction that there exists $D\in |K_{\Si}+L_2|$ and consider the
image $C$ of $D$ in $\pp^2$; $C$ is a cubic containing $P_1,\dots,P_6$ which
has a double point at $P_5$ and $P_6$. By Bezout's theorem, $\De_3$ is
contained in $C$ and thus $C=\De_3+Q$, where $Q$ is a conic containing
$P_1,\dots,P_6$, which is impossible. By similar (easier) arguments, one
shows that $h^0(\Si,K_{\Si}+L_1)=h^0(\Si, K_{\Si}+L_3)=0$, and thus
$p_g(S)=p_g(X)=0$. By the projection formula for a finite flat morphism,
 
\begin{multline*}
H^0(X,2K_X)= \\
H^0(\Si, {-}K_{\Si}+f_1+\sum S_j)\oplus\Big(\bigoplus_i H^0(\Si,
{-}K_{\Si}+f_1+\sum S_j-L_i)\Big),
 \end{multline*}
and $\Ga$ acts on $H^0(\Si, {-}K_{\Si}+f_1+\sum S_j-L_i)$ via the character
$\chi_i$. We have $h^0(\Si,{-}K_{\Si}+f_1+\sum S_j)=h^0(\Si,{-}K_{\Si}+f_1)$,
since
 
$$ S_j({-}K_{\Si}+f_1+\sum S_i)=-2 \quad\text{for $i=1,\dots,4$;}
 $$
in addition, $h^0(\Si,{-}K_{\Si}+f_1)=7$, since $\Si$ is rational,
$2f_1+f_2+f_3$ has arithmetic genus 7, and
${-}K_{\Si}+f_1=K_{\Si}+2f_1+f_2+f_3$. Since $p_2(S)=8$, there is a value
$i\in\{1,2,3\}$ such that $h^0(\Si,{-}K_{\Si}+\sum S_j+f_1-L_i)=1$ and
$h^0(\Si,{-}K_{\Si}+\sum S_j+f_1-L_k)=0$ for $k\ne i$.
Actually, an argument similar to that used for computing $p_g(S)$ shows
that
 \begin{align*}
 h^0({-}K_{\Si}+\sum S_j+f_1-L_1) &=h^0(\sum S_j+e_4)=1, \\
 h^0({-}K_{\Si}+\sum S_j+f_1-L_2) &=
 h^0(3l-e_1-2e_2-e_3-2e_4-e_5-e_6)=0, \\
 h^0({-}K_{\Si}+\sum S_j+f_1-L_3) &=
 h^0(5l-e_1-2e_2-e_3-3e_4-3e_5-3e_6)=0.
 \end{align*}
It follows that the bicanonical map $\fie\colon S\to \pp^7$ is composed with
the involution $\ga_1$ but not with $\ga_2$ and $\ga_3$. Since
$|2K_S|\supset \pi^*|{-}K_{\Si}|$ and the map $\Si \to\pp^3$ induced
by $|{-}K_{\Si}|$ is birational, it follows that $\fie$ has degree $2$.

\medskip

The remaining examples are obtained using the following construction due
to Beauville (see \cite [p.~123, Ex.~4] {bv} and cf.\ \cite{dolg}). Let
$C_1,C_2$ be curves of genus $g_1,g_2$, and assume that a group $G$ of
order $(g_1-1)(g_2-1)$ acts on $C_1,C_2$ so that $C_i/G$ is isomorphic to
$\pp^1$ for $i=1,2$; write $p_i\colon C_i\to\pp^1$ for the projections onto
the quotients and $p\colon C_1\times C_2 \to\pp^1\times\pp^1$ for the
product of $p_1$ and $p_2$. Thus $p$ is a Galois cover with group $G\times G$.
Assume in addition that there exists an automorphism $\psi\in\Aut G$ whose
graph $\Ga=\Ga_{\psi}\subset G\times G$ acts freely on $C_1\times C_2$. Then
set $S= (C_1\times C_2)/\Ga$ and denote by $q\colon C_1\times C_2\to S$ the
quotient map and by $\pi\colon S\to \pp^1\times\pp^1$ the map induced by $p$.
If $G$ is Abelian, then $\pi$ is a $G$-cover. The surface $S$ is minimal and
of general type since $C_1\times C_2$ is minimal of general type and $q$ is
\'etale. Since $\Ga$ acts freely,
$\chi(\Oh_{C_1\times C_2})=|G|\chi(\Oh_S)$ and $K^2_{C_1\times C_2}=|G|K_S^2$,
namely $\chi(\Oh_S)=1$, $K^2_S=8$. The irregularity $q(S)$ equals the
dimension of the $\Ga$-invariant subspace of
$H^0(C_1\times C_2, \Omega^1_{C_1\times C_2})\iso
H^0(C_1,\omega_{C_1})\oplus H^0(C_2,\omega_{C_2})$. Since $C_1/G$ and $C_2/G$
are both rational and $\psi$ is an automorphism, it follows that $q(S)=0$,
and thus $p_g(S)=0$.
\end{ex}

\begin{ex}\label{ex!2} As far as we know, this is a new example. In this
case $G=\Z_2^3$, $g_1=5$, $g_2=3$. We denote by $\ga_1,\ga_2,\ga_3$ the
standard generators of $G$ and by $\chi_1,\chi_2,\chi_3$ the dual basis of
the group of characters $G^*$. To construct the $G$-cover
$p_i\colon C_i\to\pp^1$ we have to specify (cf.\ \cite[Propositions~2.1
and~3.1]{ritaabel}):

 \begin{enumerate}
 \item a divisor $D_{\ga}$ on $\pp^1$ for each nonzero $\ga\in G$;
 \item line bundles $L_1,L_2,L_3$ on $\pp^1$ satisfying
 \[
 2L_i\equiv \sum_{\ga} \epsi_i(\ga)D_{\ga},
 \quad\text{where}\quad
 \begin{cases}
 \epsi_i(\ga)=0 & \text{if } \ga\in \ker\chi_i, \\
\epsi_i(\ga)=1 & \text{otherwise.}
 \end{cases}
 \]
 \end{enumerate}

To construct $p_1\colon C_1\to \pp^1$, we choose distinct points
$P_1,\dots,P_6\in\pp^1$ and set $D_{\ga_1}=P_1+P_2$, $D_{\ga_2}=P_3+P_4$,
$D_{\ga_3}=P_5+P_6$, $D_{\ga}=0$ for $\ga\ne\ga_i$, and
$L_1=L_2=L_3=\Oh_{\pp^1}(1)$. The curve $C_1$ is smooth connected of genus~5.
To construct $p_2\colon C_2\to \pp^1$, we choose distinct points
$Q_1,\dots,Q_5\in\pp^1$ and set $D_{\ga_1}=Q_1$, $D_{\ga_2}=Q_2$,
$D_{\ga_1+\ga_2}=Q_3$, $D_{\ga_3}=Q_4+Q_5$, $D_{\ga}=0$ for the remaining
nonzero elements of $G$, and $L_1=L_2=L_3=\Oh_{\pp^1}(1)$. The curve $C_2$ is
smooth connected of genus $3$. Define $\psi\in\Aut G$ by
 \[\ga_1\mapsto
 \ga_1+\ga_3,\quad
 \ga_2\mapsto \ga_2+\ga_3,\quad
 \ga_3\mapsto\ga_1+\ga_2+\ga_3.
 \]
In the above notation, $\pi\colon S\to\pp^1\times\pp^1$ is a $G$-cover and
$2K_S=\pi^*\Oh_{\pp^1\times\pp^1}(2,1)$. By the projection formula, we have
 \[
 H^0(S,2K_S)=\bigoplus_{\chi\in\Ga^{\perp}}
 H^0(\pp^1\times\pp^1, \Oh_{\pp^1\times\pp^1}(2,1)\tensor M_{\chi}\inv),
 \]
where $M_{\chi}\inv$ is the eigensheaf of $\pi_*\Oh_S$ corresponding to
$\chi\in \Ga^{\perp}\iso G^*$, and $(G\times G)/\Ga\iso G$ acts on
$H^0(\pp^1\times\pp^1, \Oh_{\pp^1\times\pp^1}(2,1)\tensor M_{\chi}\inv)$ via
$\chi$. The $M_{\chi}$ are line bundles on $\pp^1\times\pp^1$ that can be
determined using \cite[(2.15)]{ritaabel}. One checks that
$H^0(\pp^1\times\pp^1,\Oh_{\pp^1\times\pp^1}(2,1)\tensor M_{\chi}\inv)$ is
nonzero only for the elements of $\Ga^{\perp}$ that are orthogonal to
$(0,\ga_3)\in G\times G$. It follows that the bicanonical map of $S$ is
composed with the involution induced on $S$ by $(0,\ga_3)$, and thus it has
degree~2 by Theorem~\ref{main}.
\end{ex}

\begin{ex}\label{ex!3} As for Example \ref {ex!1}, this is due to Inoue
\cite[p.~317]{inoue}, and arises as the quotient of a complete intersection
in the product of 4 elliptic curves by a free group action. Here we give a
construction in the style of Beauville as explained above which is more
suitable for our purpose. Let $\ga_1,\dots,\ga_4$ be a basis of $G=\Z_2^4$ and
$\chi_1,\dots,\chi_4$ the dual basis of $G^*$; set
$\ga_0=\ga_1+\ga_2+\ga_3+\ga_4$. We construct $C_i$ as $G$-covers of
$\pp^1$ for $i=1,2$. As in Example \ref {ex!2}, for this, we specify (cf.\
\cite[Propositions~2.1 and~3.1]{ritaabel}):
 \begin{enumerate}
 \item a divisor $D_{\ga}$ of $\pp^1$ for every nonzero $\ga\in G$;
 \item line bundles $L_1,\dots,L_4$ of $\pp^1$ satisfying
 \[
 2L_i\equiv \sum_{\ga} \epsi_i(\ga)D_{\ga}, \quad\text{where}\quad
 \begin{cases}
 \epsi_i(\ga)=0 & \text{if }\ga\in \ker\chi_i \\
 \epsi_i(\ga)=1 & \text{otherwise.}
 \end{cases}
 \]
 \end{enumerate}

Choose distinct points $P_0,\dots,P_4\in\pp^1$ and set $D_{\ga_i}=P_i$ for
$i=0,\dots,4$, $D_\ga=0$ if $\ga\ne \ga_i$, and $L_i=\Oh_{\pp^1}(1)$ for
$i=1,\dots,4$. We write $p_1\colon C_1\to\pp^1$ for the corresponding
$G$-cover. Then $C_1$ is a smooth connected curve of genus $5$. We construct
the curve $C_2$ in the same way, starting from points $Q_0,\dots,Q_4\in\pp^1$.

Let $\psi\in\Aut G$ be the automorphism:
 \[
 \ga_1\mapsto \ga_1+\ga_3,\quad \ga_2\mapsto \ga_2+\ga_4,\quad
 \ga_3\mapsto \ga_1+\ga_4,\quad \ga_4\mapsto \ga_1+\ga_3+\ga_4.
 \]
In the above notation, $\pi\colon S\to\pp^1\times\pp^1$ is a $G$-cover
and $2K_S=\pi^*\Oh_{\pp^1\times\pp^1}(1,1)$.  Since
$\Oh_{\pp^1\times\pp^1}(1,1)$ is  very  ample, it follows that the
bicanonical map $\fie$ of $S$ is  birational if and only if it is not
composed with an involution
$\ga$ of the Galois group
$G$ of
$\pi$.  To check that this is indeed
the case, we use the projection formula
 \[
 H^0(S,2K_S)=\bigoplus_{\chi\in \Ga^{\perp}}H^0(\pp^1\times\pp^1,
\Oh_{\pp^1\times\pp^1}(1,1)\tensor M_{\chi}\inv),
 \]
where $M_{\chi}\inv$ is the eigensheaf of $\pi_*\Oh_S$ corresponding to
$\chi\in \Ga^{\perp}\iso G^*$, and $(G\times G)/\Ga\iso G$ acts on
$H^0(\pp^1\times\pp^1,\Oh_{\pp^1\times\pp^1}(1,1)\tensor M_{\chi}\inv)$
via $\chi$. The $M_{\chi}$ are line bundles on $\pp^1\times\pp^1$ that can
be determined using \cite[(2.15)]{ritaabel}. It turns out that no $\ga\in G
\setminus\{0\}$ acts trivially on $H^0(S, 2K_S)$ and thus $\fie$ is
birational.
\end{ex}

\begin{ex}\label{ex!4} This example is due to Beauville and appears in
\cite[p.~123, Ex.~4]{bv}, where the group action is not described
explicitly, and in \cite{dolg}.  The assertion concerning the group action in  \cite{dolg} is not correct,
since the group action described  is not free.

In this case $g=6$, $G=\Z_5^2$, and $C_1=C_2=\{x^5+y^5+z^5=0\}\subset\pp^2$
is the Fermat quintic. If $\epsi$ is a primitive 5th root of 1, then
$(1,0)\in G$ acts on $C$ by $(x:y:z)\mapsto(\epsi x:y:z)$ and $(0,1)$ acts
by $(x:y:z)\mapsto(x:\epsi y:z)$. Let $\psi$ be the automorphism of $G$
taking $(1,0)\mapsto(1,-1)$ and $(0,1)\mapsto(1,2)$.

We compute the degree of the bicanonical map of $S$ by writing down an
explicit basis of the $\Ga$-invariant subspace of
$H^0(C\times C,2K_{C\times C})$. Take homogeneous coordinates
$(x:y:z;\ x_1:y_1:z_1)$ on $\pp^2\times \pp^2\supset C\times C$; using the
fact that a regular $1$-form on $C$ is the residue of a rational form
$\frac{g(x,y,z)}{x^5+y^5+z^5}\,dx\wedge dy \wedge dz$ for $g$ homogeneous
of degree~2, we see that $(a,b)\in G$ acts on bicanonical forms on
$C\times C$ by:
 \begin{multline*}
 x^iy^jz^{4-i-j}x_1^{\alpha}y_1^{\beta}z_1^{4-\alpha-\beta} \mapsto
 \epsi^l x^iy^jz^{4-i-j}x_1^{\alpha}y_1^{\beta}z_1^{4-\alpha-\beta},\\
 \text{where} \quad l=a(2+i+\alpha-\beta)+b(3+j+\alpha+2\beta)
 \end{multline*}
Thus the following is a basis of $H^0(Y, 2K_Y)^{\mathrm{inv}}$:
 \begin{gather*}
 x^4y_1z_1^3, \quad y^3zy_1^2z_1^2, \quad xyz^2y_1^3z_1, \quad
 x^2yzx_1z_1^3, \quad z^4x_1y_1^3,\\
 xz^3x_1^2z_1^2, \quad x^3yx_1^2y_1^2, \quad y^4x_1^3z_1, \quad
 xy^2zx_1^3y_1.
 \end{gather*}
The subfield of $\C(S)$ generated by  ratios of these monomials is  the
function field $\C(\Si)$ of the bicanonical image $\Si$ of $S$. 
The map $\pi\colon S\to\pp^1\times\pp^1$
identifies   $\C(\pp^1\times\pp^1)$ with the subfield of $\C(S)$
generated by
$x^5z^{-5}$ and $x_1^5z_1^{-5}$. The extension $\C(S)\supset \C(\pp^1\times\pp^1)$ is
Galois with Galois group $G=\Z_5^2$. 
We observe:
\[
x^5z^{-5}=( x^3yx_1^2y_1^2)(x^4y_1z_1^3)(x^2yzx_1z_1^3)\inv(z^4x_1y_1^3)\inv
 \]
and
 \[
x_1^5z_1^{-5}=(z^4x_1y_1^3)^2(y^3zy_1^2z_1^2)(x^3yx_1^2y_1^2)^2(x^2yzx_1z_1^3)
\inv(xyz^2y_1^3z_1)^{-4}.
 \]
It follows that $\C(\Si)\supset \C(\pp^1\times\pp^1)$. Now one checks that no
element of the Galois group $G=\Z_5^2$ of $\C(S)$ over  $\C(\pp^1\times\pp^1)$ acts
trivially on $\C(\Si)$. It follows that $\C(\Si)=\C(S)$, namely $\fie$ is birational.
\end{ex}
\paragraph{Acknowledgements} We are  grateful to the
referee for his help
in improving the presentation and for suggesting the 
calculation in Example \ref {ex!4}.\par
This work began during the final meeting in
Torino of the EU HCM network AGE contract ERBCHRXT 940557.
The research has been partly supported by the Italian P.I.N.
1997 Geo\-metria algebrica, algebra commutativa e aspetti
computazionali and by Project PRAXIS XXI 2/2.1/MAT/73/94.
The first author is a member of CMAF and of the Departamento
de Matem\'atica da Faculdade de Ci\^encias da Universidade
de Lisboa and the second author is a member of GNSAGA of CNR.

\bigskip

\begin{tabbing} 1749-016 Lisboa, PORTUGALxxxxxxxxx\= 56127 Pisa,
ITALY \kill Margarida Mendes Lopes \> Rita Pardini\\
CMAF \> Dipartimento di Matematica\\
 Universidade de Lisboa \> Universit\a`a di Pisa \\ Av. Prof. Gama Pinto, 2 \>
Via
Buonarroti 2\\ 1649-003 Lisboa, PORTUGAL \> 56127 Pisa, ITALY\\
mmlopes@lmc.fc.ul.pt \> pardini@dm.unipi.it
\end{tabbing}

\end{document}